\documentclass[10pt]{article}
\usepackage{color}
\usepackage{amssymb,amsmath,graphicx,acronym,amsfonts}
\usepackage[square, comma, sort&compress, numbers]{natbib}

\newtheorem{corollary}{Corollary}[section]
\newtheorem{theorem}{Theorem}[section]
\newtheorem{lemma}{Lemma}[section]
\newtheorem{definition}{Definition}[section]
\newtheorem{proposition}{Proposition}[section]
\newtheorem{example}{Example}[section]
\newtheorem{assum}{Assumption}[section]
\newtheorem{algo}{Algorithm}[section]
\newtheorem{Remark}{Remark}[section]

\def\bc{\begin{corl}}
\def\bc{\end{corl}}
\def\ba{\begin{algo}}
\def\ea{\end{algo}}
\def\br{\begin{Remark}}
\def\er{\end{Remark}}
\def\bs{\begin{assum}}
\def\es{\end{assum}}
\def\bt{\begin{theorem}}
\def\et{\end{theorem}\vskip 3pt}
\def\bl{\begin{lemma}}
\def\el{\end{lemma}}
\def\ep{\end{proposition}}
\def\bp{\begin{proposition}}
\def\qed{\hfill{$\Box$}\vskip 5pt}
\def\be{\begin{example}}
\def\ee{\end{example}}
\def\bd{\begin{definition}}
\def\ed{\end{definition}}
\def\bc{\begin{corollary}}
\def\ec{\end{corollary}}
\def\proof{\noindent\it Proof. \hspace{1mm}\rm}
\input epsf
\begin{document}
\title{\bf \Large Copositive Tensor Detection and Its Applications in Physics and Hypergraphs}
\author{Haibin Chen\thanks{School of Management Science, Qufu Normal University, Rizhao, Shandong, P.R. China. Email: chenhaibin508@163.com. This author's work was supported by the National Natural Science Foundation of China (Grant No. 11601261), and Natural Science Foundation of Shandong Province (Grant No. ZR2016AQ12). } \quad
Zheng-Hai Huang
\thanks{Department of Mathematics, School of Science, Tianjin University, Tianjin 300072, P.R. China. 
Email: huangzhenghai@tju.edu.cn. This author's work was supported by the National Natural Science Foundation of China (Grant No. 11431002).} \quad
Liqun Qi
\thanks{Department of Applied Mathematics, The Hong Kong Polytechnic University, Hung Hom,
Kowloon, Hong Kong. Email: maqilq@polyu.edu.hk. This author's work
was supported by the Hong Kong Research Grant Council (Grant No.
PolyU 501212, 501913, 15302114 and 15300715).}}
\date{}
\maketitle
\vspace{-0.6cm}
\begin{abstract}
Copositivity of tensors plays an important role in vacuum stability of a general scalar potential, polynomial optimization, tensor complementarity problem and
tensor generalized eigenvalue complementarity problem. In this paper, we propose a new algorithm for testing copositivity of high order tensors, and then present applications of the algorithm in physics and hypergraphs. For this purpose, we first give several new conditions for copositivity of tensors based on the representative matrix of a simplex. Then a new algorithm is proposed with the help of a proper convex subcone
of the copositive tensor cone, which is defined via the copositivity of Z-tensors. Furthermore, by considering a sum-of-squares program problem, we define two new subsets of the copositive tensor cone and discuss their convexity. As an application of the proposed algorithm,
we prove that the coclique number of a uniform hypergraph is equivalent with an optimization problem over the completely positive tensor cone,
which implies that the proposed algorithm can be applied to compute an upper bound of the coclique number of a uniform hypergraph.  Then we study another application of the proposed algorithm on particle physics in testing copositivity of some potential fields.
At last, various numerical examples are given to show the performance of the algorithm.

\medskip

\noindent{\bf Keywords:} Symmetric tensor; Strictly copositive; Positive semi-definiteness; Simplicial partition; Particle physics; Hypergraphs

\noindent{\bf AMS Subject Classification(2010):} 65H17, 15A18, 90C30.

\end{abstract}

\newpage
\section{Introduction}
Copositivity of high order tensors has received a growing amount of interest in
vacuum stability of a general scalar potential \cite{KK16}, tensor complementarity problem \cite{CQW, SQ15, SQ,BHW16,WHB16},
tensor generalized eigenvalue complementarity problem \cite{Ling16} and polynomial optimization problems \cite{Pena14,song2016}.
A symmetric tensor is called copositive if it generates a multivariate form taking nonnegative values over the nonnegative orthant \cite{qlq2013}.
In the literature, copositive tensors constitute a large class of tensors
that contain nonnegative tensors and several kinds of structured tensors in the even order symmetric case such as $M$-tensors, diagonally dominant tensors and so on \cite{CCLQ, chen14,Ding15,Kannan15,LCQL,LCL,LWZ,Qi14,qi14,Zhang12}.

Recently, Kannike \cite{KK16} studied the vacuum stability of a general scalar potential of a few fields. With the help of copositive tensors and its relationship to orbit space variables, Kannike showed that how to find positivity conditions for more complicated potentials.
Then, he discussed the vacuum stability conditions of the general potential of two real scalars, without and with the Higgs boson
included in the potential \cite{KK16}. Furthermore, explicit vacuum stability conditions for the two Higgs doublet model were given, and a short overview
of positivity conditions for tensors of quadratic couplings were established via tensor eigenvalues.

In \cite{Pena14}, Pena et al. provided a general characterization of polynomial optimization problems that can be formulated as a conic program over the cone of completely positive tensors. It is known that the cone of completely positive tensors has a natural associated dual cone of copositive
tensors \cite{QXX14}. In light of this relationship, any completely positive program stated in \cite{Pena14} has a natural dual conic program over the cone of copositive tensors. As a consequence of this characterization, it follows that recent related results for quadratic problems can be further strengthened and generalized to higher order polynomial optimization problems.   For completely positive tensors and their applications, also see \cite{K2015, LQ, QXX14}.  In \cite{song2016}, Song and Qi gave the concepts of Pareto $H$-eigenvalue (Pareto $Z$-eigenvalue) for symmetric tensors and proved that the minimum Pareto $H$-eigenvalue (Pareto $Z$-eigenvalue) is equivalent to the optimal value of a polynomial optimization problem. It is proved that a symmetric tensor $\mathcal{A}$ is strictly copositive if and only if every
Pareto $H$-eigenvalue ($Z$-eigenvalue) of $\mathcal{A}$ is positive, and $\mathcal{A}$ is copositive if and only if every
Pareto $H$-eigenvalue ($Z$-eigenvalue) of $\mathcal{A}$ is nonnegative \cite{song2016}. Unfortunately, it is NP-hard to compute the minimum Pareto H-eigenvalue or Pareto Z-eigenvalue of a general symmetric tensor.

On the other hand, Che, Qi and Wei \cite{CQW} showed that the tensor complementarity problem with a strictly copositive tensor has a nonempty and compact solution set. Song and Qi \cite{SQ15} proved that a real symmetric tensor is a (strictly) semi-positive  if and only if it is (strictly) copositive.  Song and Qi \cite{SQ15,SQ} obtained several results for the tensor complementarity problem with a (strictly) semi-positive tensor. Huang and Qi \cite{HQ2016} formulated an $n$-person noncooperative game as a tensor complementarity problem with the involved tensor being nonnegative. Thus, copositive tensors play an important role in the tensor complementarity problem. Besides, Ling et al. \cite{Ling16} gave an affirmative result
that the tensor generalized eigenvalue complementarity problem is solvable and has at least one solution under assumptions that the related tensor is
strictly copositive.

Thus, a challenging problem is how to check the copositivity of a given symmetric tensor efficiently?

Several sufficient conditions or necessary and sufficient conditions for copositive tensors have been presented in \cite{qlq2013,Song15}. However, it is hard to verify numerically whether a tensor is copositive or not from these conditions.
Actually, the problem to judge whether a symmetric tensor is copositive or not is NP-complete, even for the matrix case \cite{Dickinson14,Murty87}.
Very recently \cite{chen16}, Chen et al. gave some theoretical studies on various conditions for (strictly) copositive tensors; and based on some of the theoretical findings, several new criteria for copositive tensors are proposed based on the
representation of the multivariate form in barycentric coordinates with respect to the standard simplex
and simplicial partitions. It is verified that, as the partition gets finer and finer, the concerned
conditions eventually capture all strictly copositive tensors. It should be pointed out that the algorithm investigated in \cite{chen16} can be viewed as an extension of some branch-and-bound type algorithms for testing copositivity of symmetric matrices \cite{Bundfuss08,Sponsel12,Xu11}.

In this paper, with the help of sum-of-square polynomial technique, an alternative numerical algorithm for copositivity of tensors is proposed, which is established via a kinds of structured tensors and on the choice of a suitable convex subcone of copositive tensor cone.
It is proved that the coclique number of a uniform hypergraph is equivalent with an optimization problem over the completely positive tensor cone,
which is the dual cone of copositive tensors. Using this, the proposed algorithm can be applied to
compute the upper bound of the coclique number of a uniform hypergraph. Furthermore, the proposed algorithm is applied to test the copositivity of some potential fields on particle physics. The rest of this paper is organized as follows.

In Section 2, we recall some notions and basic facts about tensors and the corresponding homogeneous polynomials.
In Section 3, based on the corresponding matrix of a simplex, several new criteria
for (strictly) copositive tensors based on the simplicial subdivision are presented.
In Section 4, we propose the main numerical detection algorithm for copositive tensors based on a subcone of the copositive tensor cone.
The relationship between the iteration number of the algorithm and the number of all sub-simplex is established.
Furthermore, different candidates for the subcone are discussed. An upper bound for the coclique number of a uniform hypergraph is
given in Section 5. In Section 6, some numerical results are reported to verify the performance of the algorithms.
In Section 7, a particle physical example on vacuum stability is presented, and its copositivity of coupling tensors is tested by the proposed algorithm.
Some final remarks are given in Section 8.

To move on, we briefly mention the notation that will be used in the sequel. Let $\mathbb{R}^n$
be the $n$ dimensional real Euclidean space and and the set of all nonnegative vectors in $\mathbb{R}^n$ be denoted by $\mathbb{R}^n_+$.
The set all positive integers
is denoted by $\mathbb{N}$. Suppose $m, n\in \mathbb{N}$ are two natural numbers. Denote $[n]=\{1,2,\cdots,n\}$. Vectors are denoted by bold lowercase letters i.e. ${\bf x},~ {\bf y},\cdots$, matrices are denoted by capital letters i.e. $A, B, \cdots$, and tensors are written as calligraphic capitals such as
$\mathcal{A}, \mathcal{T}, \cdots.$ The $i$-th unit coordinate vector in $\mathbb{R}^n$ is denoted by ${\bf e_i}$. All one tensor and all one vector are
denoted by $\mathcal{E}$ and ${\bf e}$ respectively.
If the symbol $|\cdot|$ is used on a tensor $\mathcal{A}=(a_{i_1 \cdots i_m})_{1\leq i_j\leq n}$, $j=1,\cdots,m$, it denotes another
tensor $|\mathcal{A}|=(|a_{i_1 \cdots i_m}|)_{1\leq i_j\leq n}$, $j\in [m]$. If $\mathcal{B}=(b_{i_1 \cdots i_m})_{1\leq i_j\leq n}$, $j\in [m]$
is another tensor, then $\mathcal{A}\leq \mathcal{B}$ means $a_{i_1 \cdots i_m} \leq b_{i_1 \cdots i_m}$
for all $i_1,\cdots,i_m \in [n]$.
\setcounter{equation}{0}
\section{Preliminaries}

A real $m$-th order $n$-dimensional tensor $\mathcal{A}=(a_{i_1i_2\cdots i_m})$ is a multi-array of real entries $a_{i_1i_2\cdots i_m}$, where $i_j \in [n]$ for $j\in [m]$. In this paper, we always assume that $m\geq 3$ and $n\geq 2$. A tensor is said to be nonnegative if all its entries are nonnegative. If the entries $a_{i_1i_2\cdots i_m}$ are invariant under any permutation of their indices, then tensor $\mathcal{A}$ is called a symmetric tensor. In this paper, we always consider real symmetric tensors. The identity tensor $\mathcal{I}$ with order $m$ and dimension $n$ is given by $\mathcal{I}_{i_1\cdots i_m}=1$ if $i_1=\cdots=i_m$ and $\mathcal{I}_{i_1\cdots i_m}=0$ otherwise. All one tensor $\mathcal{E}$ (all one vector ${\bf e}$) is a tensor (vector) with all entries
equal one.

We denote
$$
\mathbb{S}_{m,n}:=\{\mathcal{A}: \mathcal{A} \mbox{ is an } m\mbox{-th~order }  n\mbox{-dimensional} \mbox{  symmetric tensor}\}.
$$
Clearly, $\mathbb{S}_{m,n}$ is a vector space under the addition and multiplication defined as below: for any $t \in \mathbb{R}$,
$\mathcal{A}=(a_{i_1 \cdots i_m})_{1 \le i_1,\cdots,i_m
\le n}$ and $\mathcal{B}=(b_{i_1 \cdots i_m})_{1 \le
i_1,\cdots,i_m \le n},$
\[
\mathcal{A}+\mathcal{B}=(a_{i_1 \cdots i_m}+b_{i_1 \cdots i_m})_{1
\le i_1,\cdots,i_m \le n}\quad \mbox{\rm and }\quad t
\mathcal{A}=(ta_{i_1 \cdots i_m})_{1 \le i_1,\cdots,i_m
\le n}.
\]
In addition, there are some more tensor cones that will be used in the following analysis such as copositive tensor cone ($\mathbb{COP}_{m,n}$),
completely positive tensor cone ($\mathbb{CP}_{m,n}$), nonnegative tensor cone ($\mathbb{N}^+_{m,n}$) and positive semi-definite tensor cone ($\mathbb{PSD}$).

For any $\mathcal{A}, \mathcal{B} \in \mathbb{S}_{m,n}$, we define the inner product by
$\langle \mathcal{A},\mathcal{B}\rangle:=\sum_{i_1,\cdots,i_m=1}^{n}a_{i_1 \cdots i_m}b_{i_1 \cdots i_m}$, and the corresponding norm by
$$
\|\mathcal{A}\|=\left(\langle \mathcal{A},\mathcal{A}\rangle\right)^{1/2}=\left(\sum_{i_1,\cdots,i_m=1}^{n}(a_{i_1
\cdots i_m})^2\right)^{1/2}.
$$

For any ${\bf x}\in \mathbb{R}^n$, we use $x_i$ to denote its $i$th component; and use $\|{\bf x}\|_m$ to denote the $m$-norm of ${\bf x}$.

For $m$ vectors ${\bf x},{\bf y}, \cdots, {\bf z}\in
\mathbb{R}^n$, we use ${\bf x}\circ {\bf y}\circ \cdots \circ {\bf z}$ to denote the $m$-th order $n$-dimensional symmetric rank one tensor with
\[
({\bf x}\circ {\bf y}\circ \cdots \circ {\bf z})_{i_1 i_2\cdots i_m}=x_{i_1}y_{i_2}\cdots z_{i_m}, \ \forall \, i_1,\cdots,i_m \in [n].
\]
And the inner product of a symmetric tensor and the rank one tensor is given by
$$
\langle
\mathcal{A}, {\bf x}\circ {\bf y}\circ \cdots \circ {\bf z}\rangle:=\sum_{i_1,\cdots,i_m=1}^{n}a_{i_1 \cdots
i_m}x_{i_1}y_{i_2}\cdots z_{i_m}.
$$
For $m\in \mathbb{N}$ and $k\in [m]$, we denote
$$
\mathcal{A}{\bf x}^k{\bf y}^{m-k}=\langle\mathcal{A}, \underbrace{{\bf x}\circ \cdots {\bf x}}_k\circ\underbrace{{\bf y}\circ \cdots \circ {\bf y}}_{m-k}\rangle\quad \mbox{\rm and}\quad \mathcal{A}{\bf x}^m=\langle\mathcal{A}, \underbrace{{\bf x}\circ \cdots {\bf x}}_m\rangle,
$$
then
\begin{equation}\label{e21}
\mathcal{A}{\bf x}^k{\bf y}^{m-k}=\sum_{i_1,\cdots,i_m=1}^{n}a_{i_1 \cdots i_m}x_{i_1}\cdots x_{i_k}y_{i_{k+1}}\cdots y_{i_m}
\quad \mbox{\rm and}\quad \mathcal{A}{\bf x}^m=\sum_{i_1,\cdots,i_m=1}^{n}a_{i_1 \cdots i_m}x_{i_1}\cdots x_{i_m}.
\end{equation}
For any $\mathcal{A}=(a_{i_1i_2\cdots i_m})\in \mathbb{S}_{m,n}$ and ${\bf x}\in \mathbb{R}^n$, we have $\mathcal{A}{\bf x}^{m-1}\in \mathbb{R}^n$ with
$$
(\mathcal{A}{\bf x}^{m-1})_i=\sum_{i_2,i_3,\cdots,i_m\in [n]}a_{ii_2\cdots i_m}x_{i_2}\cdots x_{i_m},~~\forall \, i\in [n].
$$

It is known that an $m$-th order $n$-dimensional symmetric tensor defines uniquely an $m$-th degree homogeneous
polynomial $f_{\mathcal{A}}({\bf x})$ on $\mathbb{R}^n$: for all ${\bf x}=(x_1,\cdots,x_n)^T
\in \mathbb{R}^n$,
$$
f_{\mathcal{A}}({\bf x})= \mathcal{A}{\bf x}^m=\sum_{i_1,i_2,\cdots, i_m\in [n]}a_{i_1i_2\cdots i_m}x_{i_1}x_{i_2}\cdots x_{i_m};
$$
and conversely, any $m$-th degree homogeneous polynomial function $f({\bf x})$ on $\mathbb{R}^n$ also corresponds uniquely a symmetric tensor. Furthermore, an even order symmetric tensor $\mathcal{A}$ is called positive semi-definite (positive definite) if $f_{\mathcal{A}}({\bf x}) \geq 0$ ($f_{\mathcal{A}}({\bf x})> 0$) for all ${\bf x}\in \mathbb{R}^n$ (${\bf x}\in \mathbb{R}^n \backslash \{\bf 0\}$).

To end this section, we introduce the notion of tensor product, which will be used in the following analysis.
\bd\label{def21}{\bf \cite{Shao13}} Let $\mathcal{A}$
($\mathcal{B}$) be an order $m\geq2$ (an order $k\geq1$) dimension $n$ tensor.
The product $\mathcal{A}\mathcal{B}$ is the following tensor $\mathcal{C}$ of order $(m-1)(k-1)+1$
with entries:
$$c_{i\alpha_1\alpha_2\cdots\alpha_{m-1}}=\sum_{i_2,\cdots,i_m\in [n_2]}a_{ii_2\cdots i_m}b_{i_2\alpha_1}\cdots b_{i_m\alpha_{m-1}},$$
where $i\in [n], \alpha_1,\alpha_2,\cdots,\alpha_{m-1} \in [n]^{k-1}$.
\ed
Here, when tensor $\mathcal{B}$ reduces to a 1st order tensor, i.e., vector of $\mathbb{R}^n$, the production
$\mathcal{A}{\bf x}$ coincides with the notation $\mathcal{A}{\bf x}^m$

\setcounter{equation}{0}
\section{Several conditions of copositive tensors}
In this section, we will give several new sufficient conditions or necessary conditions with the help of a proper subcone of
the cone of copositive tensors.
As a simplex $S$ is determined by its vertices, it also can be represented by a matrix $V_S$ whose
columns are these vertices. $V_S$ is nonsingular and unique up to a permutation of its columns.
In the following analysis, we always assume that $\mathbb{M}$ is a subcone of the copositive tensor cone $\mathbb{COP}_{m,n}$.

Before we move on, we first list the definition of copositive tensors and some notions about simplex.

\bd\label{def31} {\bf \cite{qlq2013}} Let $\mathcal{A}\in S_{m,n}$ be given. If
$\mathcal{A}{\bf x}^m\geq 0\;(\mathcal{A}{\bf x}^m>0)$ for any ${\bf x}\in \mathbb{R}^n_+\;({\bf x}\in \mathbb{R}^n_+\backslash \{{\bf0}\})$, then
$\mathcal{A}$ is called a copositive (strictly copositive) tensor.
\ed

The standard simplex with vertices ${\bf e}_1, {\bf e}_2,\cdots,{\bf e}_n$ is denoted by $S_0=\{{\bf x}\in \mathbb{R}^n_+~|~\|{\bf x}\|_1=1 \}$.
Let $S,S_1,S_2,\cdots,S_r$ be finite simplices in $\mathbb{R}^n$. The set $\tilde{S}=\{S_1,S_2,\cdots,S_r\}$ is called a simplicial partition of $S$ if it satisfies that
$$
S=\bigcup_{i=1}^rS_i\quad \mbox{\rm and}\quad \mbox{\rm int} S_i\bigcap \mbox{\rm int} S_j=\emptyset \;\, \mbox{\rm for any}\;\, i,j\in [r]\;\, \mbox{\rm with}\; i\neq j,
$$
where $\mbox{\rm int}S_i$ denotes the interior of $S_i$ for any $i\in [r]$. Let $d(\widetilde{S})$ denote the maximum diameter of a simplex in $\widetilde{S}$, which is given by
$$
d(\widetilde{S})=\max_{k\in [r]}\max_{i,j\in [n]}\|{\bf u}^k_i-{\bf u}^k_j\|_2.
$$

\bt\label{them31}
Suppose $\mathcal{A}\in \mathbb{S}_{m,n}$. Let $S_1=conv\{ {\bf u}_1,{\bf u}_2,\cdots,{\bf u}_n\}$ be a simplex, where ${\bf u}_i\in \mathbb{R}^n, i\in [n]$. Let $V=({\bf u}_1,{\bf u}_2,\cdots,{\bf u}_n)$ be the square matrix corresponding to $S_1$. If $V^T\mathcal{A}V\in \mathbb{M}$, then $\mathcal{A}{\bf x}^m\geq 0$ for all ${\bf x}\in S_1$.\et
\proof It is apparent that $V^T\mathcal{A}V$ is an $m$-th order $n$ dimensional tensor with entries
$$
\begin{aligned} (V^T\mathcal{A}V)_{i_1i_2\cdots i_m}=&\sum_{j_1,j_2,\cdots,j_m\in [n]}(V^T)_{i_1j_1}a_{j_1j_2\cdots j_m}V_{j_2i_2}
\cdots V_{j_mi_m} \\
=&\sum_{j_1,j_2,\cdots,j_m\in [n]}V_{j_1i_1}a_{j_1j_2\cdots j_m}V_{j_2i_2}
\cdots V_{j_mi_m}\\
=&\sum_{j_1,j_2,\cdots,j_m\in [n]}a_{j_1j_2\cdots j_m} ({\bf u}_{i_1})_{j_1}({\bf u}_{i_2})_{j_2}\cdots ({\bf u}_{i_m})_{j_m}\\
=&\langle\mathcal{A}, {\bf u}_{i_1} \circ{\bf u}_{i_2}\circ \cdots \circ {\bf u}_{i_m} \rangle,
\end{aligned}
$$
for all $i_1, i_2, \cdots, i_m\in [n]$.
For any ${\bf x}\in S_1$, for some $k\in [r]$, it follows that,
$$
{\bf x}=x_1 {\bf u}_1+x_2 {\bf u}_2+\cdots+x_n {\bf u}_n,~~\sum_{i=1}^n x_i=1,~~x_i\geq 0,~\forall\, i\in [n].
$$
Thus,
$$
\begin{aligned}
\mathcal{A}{\bf x}^m=&\langle\mathcal{A}, (x_1 {\bf u}_1+x_2 {\bf u}_2+\cdots+x_n {\bf u}_n)^m\rangle \\
=&\sum_{i_1,i_2,\cdots,i_m\in [n]}x_{i_1}x_{i_2}\cdots x_{i_m}\langle\mathcal{A}, {\bf u}_{i_1} \circ{\bf u}_{i_2}\circ \cdots \circ {\bf u}_{i_m} \rangle \\
=&\sum_{i_1,i_2,\cdots,i_m\in [n]}(V^T\mathcal{A}V)_{i_1i_2\cdots i_m}x_{i_1}x_{i_2}\cdots x_{i_m} \\
=&(V^T\mathcal{A}V)\bar{{\bf x}}^m \\
\geq & 0,
\end{aligned}
$$
since $\bar{{\bf x}}^T=(x_1,x_2,\cdots,x_n)\in \mathbb{R}^n_+$ and $V^T\mathcal{A}V$ is copositive, and hence, the desired result follows.  \qed

\bc\label{corol31}
Let $\mathcal{A}\in \mathbb{S}_{m,n}$ be given. Suppose $\widetilde{S}=\{S_1,S_2,\cdots,S_r\}$ is a simplicial partition of simplex $S_0=\{{\bf x}\in \mathbb{R}^n_+~|~\|{\bf x}\|_1=1 \}$; and the vertices of simplex $S_k$ are denoted by ${\bf u}^k_1,{\bf u}^k_2,\cdots,{\bf u}^k_n$ for any $k\in [r]$. Let $V_{S_k}=({\bf u}^k_1, {\bf u}^k_2, \cdots, {\bf u}^k_n)$ be the matrix corresponding to simplex $S_k$ for any $k\in [r]$. Then $\mathcal{A}$ is copositive if $V_{S_k}^T\mathcal{A}V_{S_k}\in \mathbb{M}$ for all $k\in [r]$. \ec

In fact, Corollary \ref{corol31} is a generalization of the sufficient condition proposed in \cite{chen16}.
In order to give a necessary condition for the strictly copositive tensor,
we cite a useful result below.

\bl\label{lema31}{\bf \cite{chen16}} Let $\mathcal{A}\in \mathbb{S}_{m,n}$ be a strictly copositive tensor. Then, there exists $\varepsilon>0$ such that for all
finite simplicial partitions $\widetilde{S}=\{S_1,S_2,\cdots,S_r\}$ of $S_0$ with $d(\widetilde{S})<\varepsilon$, it follows that
$$
\langle \mathcal{A},~{\bf u}^k_{i_1}\circ {\bf u}^k_{i_2}\circ \cdots\circ {\bf u}^k_{i_m} \rangle > 0
$$
for all $k\in [r], i_j\in [n],j\in [m]$, where ${\bf u}^k_1, {\bf u}^k_2,\cdots, {\bf u}^k_n$
are vertices of the simplex $S_k$.
\el

\bt\label{them32}Let $\mathcal{A}\in \mathbb{S}_{m,n}$ be a strictly copositive tensor. Suppose $\mathbb{M}\supseteq \mathbb{N}^+_{m,n}$. Then, there exists $\varepsilon>0$ such that for all
finite simplicial partitions $\widetilde{S}=\{S_1,S_2,\cdots,S_r\}$ of $S_0$ with $d(\widetilde{S})<\varepsilon$, it follows that
$V_{S_k}^T\mathcal{A}V_{S_k}\in \mathbb{M}$ for all $k\in [r]$, where $V_{S_k}=({\bf u}^k_1, {\bf u}^k_2,\cdots, {\bf u}^k_n)\in \mathbb{R}^{n\times n}$ and ${\bf u}^k_1, {\bf u}^k_2,\cdots, {\bf u}^k_n$ are vertices of the simplex $S_k$. \et
\proof By Lemma \ref{lema31}, we obtain that
$$
\begin{aligned} (V_{S_k}^T\mathcal{A}V_{S_k})_{i_1i_2\cdots i_m}=&\sum_{j_1,j_2,\cdots,j_m\in [n]}(V_{S_k}^T)_{i_1j_1}a_{j_1j_2\cdots j_m}(V_{S_k})_{j_2i_2}
\cdots (V_{S_k})_{j_mi_m} \\
=&\sum_{j_1,j_2,\cdots,j_m\in [n]}(V_{S_k})_{j_1i_1}a_{j_1j_2\cdots j_m}(V_{S_k})_{j_2i_2}
\cdots (V_{S_k})_{j_mi_m}\\
=&\sum_{j_1,j_2,\cdots,j_m\in [n]}a_{j_1j_2\cdots j_m} ({\bf u}^k_{i_1})_{j_1}({\bf u}^k_{i_2})_{j_2}\cdots ({\bf u}^k_{i_m})_{j_m}\\
=&\langle\mathcal{A}, {\bf u}^k_{i_1} \circ{\bf u}^k_{i_2}\circ \cdots \circ {\bf u}^k_{i_m} \rangle \\
>& 0.
\end{aligned}
$$
By Lemma \ref{lema31}, we know that $V_{S_k}^T\mathcal{A}V_{S_k}>0$ when the $\varepsilon>0$ is small enough. Thus $V_{S_k}^T\mathcal{A}V_{S_k}\in \mathbb{N}^+_{m,n}\subseteq\mathbb{M}$ for all $k\in [r]$, and the desired results hold. \qed

\bt\label{them33} Suppose $\mathcal{A}\in \mathbb{S}_{m,n}$ is copositive.
Let $S=conv\{{\bf u}_1,{\bf u}_2,\cdots,{\bf u}_n\}$ be a simplex with $\mathcal{A}{\bf u}_i^m>0$ for all $i\in [n]$.
Let $V=({\bf u}_1,{\bf u}_2,\cdots,{\bf u}_n)\in \mathbb{R}^{n\times n}$.
If there exists $\tilde{{\bf x}}\in S\backslash \{{\bf u}_1,{\bf u}_2,\cdots,{\bf u}_n\}$ such that $\mathcal{A}\tilde{{\bf x}}^m=0$,
then $V^T\mathcal{A}V$ is not strictly copositive.
\et
\proof We will prove the conclusion by contradiction. Assume $V^T\mathcal{A}V$ is strictly copositive. Then, for any ${\bf x}\in \mathbb{R}_+^n \setminus \{{\bf 0}\}$, we have
$$
\begin{aligned}
(V^T\mathcal{A}V){\bf x}^m=&\sum_{i_1,i_2,\cdots,i_m\in [n]}(V^T\mathcal{A}V)_{i_1i_2\cdots i_m}x_{i_1}x_{i_2}\cdots x_{i_m} \\
=& \sum_{i_1,i_2,\cdots,i_m\in [n]}x_{i_1}x_{i_2}\cdots x_{i_m}\langle\mathcal{A}, {\bf u}_{i_1} \circ{\bf u}_{i_2}\circ \cdots \circ {\bf u}_{i_m} \rangle \\
> & 0.
\end{aligned}
$$

Since $\tilde{{\bf x}}\neq {\bf 0}$ and ${\bf u}_1,{\bf u}_2,\cdots,{\bf u}_n$ constitute a basis of $\mathbb{R}^n$, by letting $\tilde{{\bf x}}=\tilde{{\bf x}}_1{\bf u}_1+\tilde{{\bf x}}_2{\bf u}_2+\cdots+\tilde{{\bf x}}_n{\bf u}_n$, it follows that
$$
\begin{aligned}
0=&\mathcal{A}\tilde{{\bf x}}^m \\
=&\langle \mathcal{A}, (\tilde{{\bf x}}_1{\bf u}_1+\tilde{{\bf x}}_2{\bf u}_2+\cdots+\tilde{{\bf x}}_n{\bf u}_n)^m\rangle \\
=& \sum_{i_1,i_2,\cdots,i_m\in [n]}\tilde{x}_{i_1}\tilde{x}_{i_2}\cdots \tilde{x}_{i_m}\langle\mathcal{A}, {\bf u}_{i_1} \circ{\bf u}_{i_2}\circ \cdots \circ {\bf u}_{i_m} \rangle \\
=&\sum_{i_1,i_2,\cdots,i_m\in [n]}(V^T\mathcal{A}V)_{i_1i_2\cdots i_m}\tilde{x}_{i_1}\tilde{x}_{i_2}\cdots \tilde{x}_{i_m} \\
=&(V^T\mathcal{A}V)\tilde{{\bf x}}^m,
\end{aligned}
$$
which contradicts the assumption that $V^T\mathcal{A}V$ is strictly copositive, and hence, the desired result holds. \qed

\setcounter{equation}{0}
\section{Algorithms}
The results of the preceding section naturally yield an algorithm to test whether a tensor is
copositive or not. Similar to the algorithm given in \cite{chen16}, we will present an algorithm by starting with the standard simplex in $\mathbb{R}^n_+$, and checking whether there is a vertex ${\bf v}$ with $\mathcal{A}{\bf v}^m < 0$,
or whether the copositivity criterion of Theorem \ref{them31} is satisfied.
First of all, we list the algorithm proposed in \cite{chen16} in the below, and then, discuss the relationship between iteration number and the number of all sub-simplices when the algorithm stops in finitely many iterations.
\vspace{4mm}

\begin{tabular}{@{}l@{}}
\hline
 \multicolumn{1}{c}{\bf Algorithm 1} \\
\hline
{\bf Input:} $\mathcal{A}\in \mathbb{S}_{m,n}$  \\
\qquad Set $\widetilde{S}:=\{S_1\}$, where $S_1=conv\{{\bf e}_1,{\bf e}_2,\cdots,{\bf e}_n\}$ is the standard simplex      \\
\qquad Set $k:=1$ \\
\qquad while $k\neq 0$ do         \\
\qquad\qquad set $S:=S_k=conv\{{\bf u}_1,{\bf u}_2,\cdots,{\bf u}_n\}\in \widetilde{S}$ \\
\qquad\qquad if there exists $i\in [n]$ such that $\mathcal{A}{\bf u}^m_i<0$, then     \\
\qquad\qquad\qquad return ``$\mathcal{A}$ is not copositive''       \\
\qquad\qquad else if $\langle \mathcal{A}, {\bf u}_{i_1}\circ {\bf u}_{i_2}\circ \cdots \circ {\bf u}_{i_m}\rangle \geq 0$ for all $i_1,i_2,\cdots,i_m\in [n]$, then  \\
\qquad\qquad\qquad   set $\widetilde{S}:=\widetilde{S} \backslash \{S_k\}$ and $k:=k-1$    \\
\qquad\qquad else \\
\qquad\qquad\qquad set  \\
\qquad\qquad\qquad\qquad $S_k:=conv\{{\bf u}_1,\cdots,{\bf u}_{p-1},{\bf v},{\bf u}_{p+1},\cdots,{\bf u}_n\}$;  \\
\qquad\qquad\qquad\qquad $S_{k+1}:=conv\{{\bf u}_1,\cdots,{\bf u}_{q-1},{\bf v},{\bf u}_{q+1},\cdots,{\bf u}_n\}$, \\
\qquad\qquad\qquad\qquad where ${\bf v}=\frac{{\bf u_p}+{\bf u_q}}{2},\, [p,q]=arg\max_{i,j\in [n]}\|{\bf u_i}-{\bf u_j}\|_2$ and $p<q$.\\
\qquad\qquad\qquad set $\widetilde{S}:=\widetilde{S}\backslash \{S\}\bigcup \{S_k, S_{k+1}\}$ and $k:=k+1$ \\
\qquad\qquad end if \\
\qquad end while \\
\qquad return `` $\mathcal{A}$ is copositive.'' \\
{\bf Output:} ``$\mathcal{A}$ is copositive'' or ``$\mathcal{A}$ is not copositive''.\\
\hline
\end{tabular}

\vspace{4mm}

For the standard simplex $S=conv\{{\bf e}_1,{\bf e}_2,\cdots,{\bf e}_n\}$ and its simplicial partition $\tilde{S}=\{S_1,S_2,\cdots,S_r\}$,
any simplex $S_i, i\in [r]$ is called a sub-simplex of $S$.
\bp For a given tensor $\mathcal{A}\in \mathbb{S}_{m,n}$, if Algorithm 1 stops in the $k$-th iteration, $k\geq 2$, then
the number of all sub-simplices need to be checked during the whole running process is $d=\frac{k+1}{2}$. Furthermore, if $\mathcal{A}$ is nonnegative or $\mathcal{A}$ has negative diagonal entries, then $k=1$.
\ep
\proof By conditions, the algorithm stops in $k$-th iteration.
Assume the original standard simplex is cut $t$ times from beginning to the end, by the fact that
the number of simplices will increase 1 when it is cut one time, so we have
$$d=t+1~~ {\bf and}~~ t+t+1=k,$$
which implies that $d=\frac{k+1}{2}$. When $\mathcal{A}$ is nonnegative or $\mathcal{A}$ has negative diagonal entries, it means that
$$
\langle \mathcal{A},~{\bf e_{i_1}}\circ {\bf e_{i_2}}\circ \cdots \circ {\bf e_{i_m}} \rangle=a_{i_1i_2\cdots i_m} \geq 0, ~\forall~ i_1,i_2,\cdots,i_m \in [n],
$$
or $\mathcal{A}{\bf e}_i^m<0$ for some $i\in [n]$. Thus, this algorithm will stop in one iteration, i.e., $k=1$, and hence, the desired result holds. \qed

We now list the main algorithm related to a convex cone $\mathbb{M}$, which is a subcone of the copositive tensor cone.
Then, different choice for the subcone $\mathbb{M}$ are discussed in detail.
\vspace{4mm}

\begin{tabular}{@{}l@{}}
\hline
 \multicolumn{1}{c}{\bf Algorithm 2} \\
\hline
{\bf Input:} $\mathcal{A}\in \mathbb{S}_{m,n}$  \\
\qquad Set $\widetilde{S}:=\{S_1\}$, where $S_1=conv\{{\bf e}_1,{\bf e}_2,\cdots,{\bf e}_n\}$ is the standard simplex      \\
\qquad Set $k:=1$ \\
\qquad while $k\neq 0$ do         \\
\qquad\qquad set $S:=S_k=conv\{{\bf u}_1,{\bf u}_2,\cdots,{\bf u}_n\}\in \widetilde{S}$ \\
\qquad\qquad let $V=({\bf u}_1,{\bf u}_2,\cdots,{\bf u}_n)$ be the square matrix corresponding to $S$ \\
\qquad\qquad if there exists $i\in [n]$ such that $\mathcal{A}{\bf u}^m_i<0$, then     \\
\qquad\qquad\qquad return ``$\mathcal{A}$ is not copositive''       \\
\qquad\qquad else if $V^T\mathcal{A}V\in \mathbb{M}$, then $\widetilde{S}=\widetilde{S} \backslash \{S_k\}$ and $k:=k-1$         \\\
\qquad\qquad else \\
\qquad\qquad\qquad set  \\
\qquad\qquad\qquad\qquad $S_k:=conv\{{\bf u}_1,\cdots,{\bf u}_{p-1},{\bf v},{\bf u}_{p+1},\cdots,{\bf u}_n\}$;  \\
\qquad\qquad\qquad\qquad $S_{k+1}:=conv\{{\bf u}_1,\cdots,{\bf u}_{q-1},{\bf v},{\bf u}_{q+1},\cdots,{\bf u}_n\}$, \\
\qquad\qquad\qquad\qquad where ${\bf v}=\frac{{\bf u_p}+{\bf u_q}}{2},\, [p,q]=arg\max_{i,j\in [n]}\|{\bf u_i}-{\bf u_j}\|_2$ and $p<q$.\\
\qquad\qquad\qquad set $\widetilde{S}:=\widetilde{S}\backslash \{S\}\bigcup \{S_k, S_{k+1}\}$ and $k:=k+1$ \\
\qquad\qquad end if \\
\qquad end while \\
\qquad return `` $\mathcal{A}$ is copositive.'' \\
{\bf Output:} ``$\mathcal{A}$ is copositive'' or ``$\mathcal{A}$ is not copositive''.\\
\hline
\end{tabular}
\vspace{4mm}

\br By the analysis in Section 4, whether or not the algorithm does terminate depends
on the input tensor $\mathcal{A}$.

{\bf (i)} If the input tensor $\mathcal{A}$ is not copositive, then the algorithm terminates. In this case, it
does not matter which set $\mathbb{M}$ is used.

{\bf (ii)} If the input tensor $\mathcal{A}$ is strictly copositive and $\mathbb{M}\supseteq \mathbb{N}_{m, n}^+$, then the algorithm
terminates in finitely many iterations.

{\bf (iii)} If the input tensor $\mathcal{A}$ is copositive but not strictly copositive, then the algorithm may or may
not terminate.
\er

An important issue which influences the number of iterations and the runtime of Algorithm 2
is the choice of the set $\mathbb{M}$. The desirable properties of the set $\mathbb{M}$ used here can be summarized such that,
for any given symmetric tensor $\mathcal{A}$, we can easily check whether $\mathcal{A}\in \mathbb{M}$; and $\mathbb{M}$ is a
subcone of the copostive tensor cone that is as large as possible.

\subsection{The choice that $\mathbb{M}=\mathbb{N}^+_{m,n}$}
The first choice one may consider easily is $\mathbb{M} = \mathbb{N}^+_{m,n}$. However, this is not always
desirable. To check whether a symmetric tensor belongs to $\mathbb{N}^+_{m,n}$ does not take much effort, but the nonnegative
tensor cone is a quite bad approximation of the copositive tensor cone. So each iteration of the algorithm
is cheap but the number of iterations may tend to be large. On the other side, it should be noted that, in Algorithm 1, the following inequality
$$\langle \mathcal{A}, {\bf u}_{i_1}\circ {\bf u}_{i_2}\circ \cdots \circ {\bf u}_{i_m}\rangle \geq 0,~\forall~i_1,i_2,\cdots,i_m\in [n],$$
exactly imply $V^T\mathcal{A}V\in \mathbb{M}=\mathbb{N}^+_{m,n}$ and the converse may not be true in general.

\subsection{An alternative choice that $\mathbb{M}$ is related with $Z$-tensors}
In order to choose a good approximation of the copositive tensor cone, we first recall the matrix case. It is obvious that $\mathbb{PSD}+\mathbb{N}^+_{m,n}$ is a good approximation of copositive matrix cone \cite{Sponsel12,Tanaka16}. The problem of testing a given matrix whether or not belongs to $\mathbb{PSD}+\mathbb{N}^+_{m,n}$ can be solved by solving the following doubly nonnegative program
$$
\begin{aligned}
{\bf Minimize}~~&~\langle A, X\rangle \\
{\bf subject~to}~&~\langle I_n, X\rangle=1,~X\in \mathbb{PSD}\cap\mathbb{N}^+_{m,n},
\end{aligned}
$$
which can be expressed as a semidefinite program. Thus, the set $\mathbb{PSD}+\mathbb{N}^+_{m,n}$ is a rather large and tractable convex subcone of $\mathbb{COP}_{m,n}$. However, solving the doubly nonnegative problem takes an awful lot of time \cite{Sponsel12,Yoshi10} and does not make for
a practical implementation. To overcome the drawback, in \cite{Sponsel12}, more easily tractable subcones of the copostive matrix cone are proposed
such that
$$
\mathbb{H}=\{A\in \mathbb{S}_n~|~A-N(A)\in \mathbb{PSD}\},
$$
where $N(A)$ is a square matrix such that
$$
N(A)_{ij}=\left\{
\begin{array}{cll}
A_{ij}~&~~A_{ij}>0~{\bf and}~i\neq j, \\
0~~&~~{\bf otherwise}.
\end{array}
\right.
$$
Here, $A-N(A)$ is a $Z$-matrix. Stimulated by this method and the notion of $Z$-tensors \cite{Zhang12}, we now extend this subcone to the high order tensor case.

As we all know that checking the positive semi-definiteness of a general symmetric tensor is NP-hard \cite{Lim2013}.
But for some symmetric tensors with special structure, it may not be NP-hard. Polynomial time algorithm for checking the positive semi-definiteness of $Z$-tensors was established in \cite{CLQ2016}. Now, for a given tensor $\mathcal{A}=(a_{i_1i_2\cdots i_m})$ with order $m$ and dimension $n$, let
$$
N(\mathcal{A})_{i_1i_2\cdots i_m}=\left\{
\begin{array}{cll}
a_{i_1i_2\cdots i_m}~ & ~~a_{i_1i_2\cdots i_m}>0~~{\bf and}~~\delta_{i_1i_2\cdots i_m}=0, \\
0~~ & ~~{\bf otherwise},
\end{array}
\right.
$$
where $\delta_{i_1i_2\cdots i_m}=1$ if and only if $i_1=i_2=\cdots=i_m$, otherwise $\delta_{i_1i_2\cdots i_m}=0$. Next, we will
consider a new subcone of copositive tensor cone from two cases.

\noindent{\bf (I)} When $m$ is even, we define the set $\mathbb{H}_1$ such that
$$
\mathbb{H}_1=\{\mathcal{A}\in \mathbb{S}_{m,n}~|~\mathcal{A}-N(\mathcal{A})\in \mathbb{PSD}\}.
$$
Here $\mathcal{A}-N(\mathcal{A})$ is an even order symmetric $Z$-tensor. For any ${\bf x}\in \mathbb{R}^n$, suppose $f_1({\bf x})=(\mathcal{A}-N(\mathcal{A})){\bf x}^m$ and the minimum
$H$-eigenvalue of $\mathcal{A}-N(\mathcal{A})$ is denoted by $\lambda_{{\rm min}}(\mathcal{A}-N(\mathcal{A}))$. Since an even order symmetric tensor is positive semi-definite if and
only if its minimum $H$-eigenvalue is nonnegative \cite{Qi05}, by Theorem 5.1 in \cite{CLQ2016}, we know that $\mathcal{A}\in \mathbb{H}_1$ if and only if
\begin{equation}\label{e41}
\lambda_{{\rm min}}(\mathcal{A}-N(\mathcal{A}))= \max_{\mu, r \in \mathbb{R}}\{\mu: f_1({\bf x})-r (\|{\bf x}\|_m^m-1)-\mu \in \Sigma^2_m[{\bf x}]\}\geq 0,
\end{equation}
where $\Sigma_m^2[{\bf x}]$ is the set of all SOS polynomials with degree at most $m$.
It is easy to know that the sums-of-squares problem (\ref{e41})
can be equivalently rewritten as a semi-definite programming problem (SDP), and so, can be solved
efficiently. Indeed, this conversion can be done by using the commonly used Matlab Toolbox YALMIP \cite{Yalmip1,Yalmip2}. The simple code using {\rm YALMIP} is appended as follows

\begin{verbatim}
sdpvar x1 x2 ... xn r mu
f=f1(x);
g = [(x1^m+x2^m+...+xn^m)-1];
F = [sos(f1-mu-r*g)];
solvesos(F,-mu,[],[r;mu]).
\end{verbatim}
{\bf (II)} When $m$ is odd, it is well known that there is not any nontrivial odd order positive semi-definite tensors. Thus, the subcone $\mathbb{H}_2$
can be changed to another subcone of $\mathbb{COP}_{m,n}$ such that
$$
\mathbb{H}_2=\{\mathcal{A}\in \mathbb{S}_{m,n}~|~\mathcal{A}-N(\mathcal{A})\in \mathbb{COP}_{m,n}\}.
$$
Now, let $\bar{\mathcal{A}}$ be a symmetric tensor with order $2m$ and dimension $n$ with entries such that
\begin{equation}\label{e42}
f_2({\bf x})=\bar{\mathcal{A}}{\bf x}^{2m}=\sum_{i_1,i_2,\cdots,i_m\in [n]}(a_{i_1i_2\cdots i_m}-N(\mathcal{A})_{i_1i_2\cdots i_m})x_{i_1}^2x_{i_2}^2\cdots x_{i_m}^2,~~\forall~{\bf x}\in \mathbb{R}^n.
\end{equation}
Since $\mathcal{A}-N(\mathcal{A})$ is a $Z$-tensor, it is not difficult to check that $\bar{\mathcal{A}}$ is a $Z$-tensor. By (\ref{e41}) and (\ref{e42}), we know that
$$
\mathcal{A}-N(\mathcal{A})\in \mathbb{COP}_{m,n}~\Leftrightarrow~\bar{\mathcal{A}} \in \mathbb{PSD},
$$
which implies that
\begin{equation}\label{e43}
\mathcal{A}-N(\mathcal{A})\in \mathbb{COP}_{m,n}~\Leftrightarrow~\max_{\mu, r \in \mathbb{R}}\{\mu: f_2({\bf x})-r (\|{\bf x}\|_{2m}^{2m}-1)-\mu \in \Sigma^2_{2m}[{\bf x}]\}\geq 0.
\end{equation}
Thus, for a given odd order symmetric tensor $\mathcal{A}$, the matlab code using {\rm YALMIP} to check whether $\mathcal{A}\in \mathbb{H}_2$
is listed below:
\begin{verbatim}
sdpvar x1 x2 ... xn r mu
f=f2(x);
g = [(x1^2m+x2^2m+...+xn^2m)-1];
F = [sos(f2-mu-r*g)];
solvesos(F,-mu,[],[r;mu]).
\end{verbatim}

To end this section, we show the convexity of $\mathbb{H}_1, \mathbb{H}_2$. Before that, we first cite a useful lemma.

\bl\label{lema41}{\bf \cite{Qi14}} Suppose $\mathcal{C}$, $\mathcal{B}$ are two nonnegative tensors with order $m$ and dimension $n$.
If it satisfies that $|\mathcal{B}|\leq \mathcal{C}$, then $\rho{(B)} \leq \rho{(C)}$, where $\rho{(B)}, \rho{(C)}$
are spectral radius of $\mathcal{C}$ and $\mathcal{B}$ respectively.
\el

\bl\label{lema42} Let $m\in \mathbb{N}$ be even number. Suppose $\mathcal{A}$, $\mathcal{B}$ are symmetric $Z$-tensors with order $m$ and
dimension $n$. If $\mathcal{A}\leq \mathcal{B}$ and $\mathcal{A}$ is positive semi-definite, then $\mathcal{B}$ is positive semi-definite.
\el
\proof Since $\mathcal{A}, \mathcal{B}$ are $Z$-tensors, we can find $t\in \mathbb{R}, t>0 $ such that
$$
\mathcal{A}=t\mathcal{I}-\mathcal{A}',~~\mathcal{B}=t\mathcal{I}-\mathcal{B}',
$$
where $\mathcal{A}', \mathcal{B}'$ are nonnegative tensors. It is easy to know that $\mathcal{A}' \geq \mathcal{B}'$ since $\mathcal{A} \leq \mathcal{B}$.
By Lemma \ref{lema41} and Corollary 3 in \cite{Qi05}, it follows that $\rho{(A')} \geq \rho{(B')}$ and
$$
\lambda_{min}(\mathcal{B})=t-\rho(\mathcal{B}')\geq t-\rho(\mathcal{A}')=\lambda_{min}(\mathcal{A}).
$$
Here, $\lambda_{min}(\mathcal{A}), \lambda_{min}(\mathcal{B})$ denote the minimum $H$-eigenvalues of $\mathcal{A}$ and $\mathcal{B}$ respectively.
By conditions that $\mathcal{A}$ is positive semi-definite, we obtain that
$$
\lambda_{min}(\mathcal{B})\geq \lambda_{min}(\mathcal{A})\geq 0,
$$
which implies that $\mathcal{B}$ is positive semi-definite, and hence, the desired results hold.
\qed
\bt\label{them41} Let $m\in \mathbb{N}$ be even. Then $\mathbb{H}_1$ is convex and it satisfies $\mathbb{N}^+_{m,n}\subseteq \mathbb{H}_1\subseteq \mathbb{COP}_{m,n}$.
\et
\proof The second statement is obvious by their definitions respectively, so we only need to prove the convexity of $\mathcal{H}_1$ i.e., $\mathcal{A}+\mathcal{B}\in \mathbb{H}_1$
for any $\mathcal{A}, \mathcal{B}\in \mathbb{H}_1$. Suppose $\mathcal{A}, \mathcal{B}\in \mathbb{H}_1$, by the definition of $\mathbb{H}_1$, we have
that
\begin{equation}\label{e44}
\mathcal{A}-N(\mathcal{A})\in \mathbb{PSD},~~\mathcal{B}-N(\mathcal{B})\in \mathbb{PSD}.
\end{equation}
By the fact that $N(\mathcal{A+B})\leq N(\mathcal{A})+N(\mathcal{B})$, we obtain that
$$
\mathcal{\mathcal{A+B}}-N(\mathcal{A+B})\geq \mathcal{A+B}-N(\mathcal{A})-N(\mathcal{B}).
$$
By (\ref{e44}) and Lemma \ref{lema42}, we know that $\mathcal{\mathcal{A+B}}-N(\mathcal{A+B})\in \mathbb{PSD}$, which implies
that $\mathcal{A+B}\in \mathbb{H}_1$. So, the desired result holds. \qed

For any $Z$-tensor $\eta \mathcal{I}-\mathcal{B}$, from Theorem 3.12 of \cite{Zhang12}, it follows that $\eta \mathcal{I}-\mathcal{B}$ is copositive if and only if $\eta\geq \rho(\mathcal{B})$. Similar to the proof of Theorem \ref{them41}, we know the following conclusion holds and the proof is omitted.
\bt\label{them42} Let $m\in \mathbb{N}$ be odd. Then $\mathbb{H}_2$ is convex and it satisfies $\mathbb{N}^+_{m,n}\subseteq \mathbb{H}_2\subseteq \mathbb{COP}_{m,n}$.
\et

By the analysis above, in Algorithm 2, we can choose $\mathbb{M}=\mathbb{H}_1$ in the even order case and $\mathbb{M}=\mathbb{H}_2$
in the odd order case.

\setcounter{equation}{0}
\section{An upper bound for the coclique number of an uniform hypergraph}
In this section, we show that computing the coclique number of a uniform hypergraph
can be reformulated as a linear program over the cone of completely positive tensors.
By the dual property of copositive tensor cone and completely positive tensor cone, we
present an upper bound for the coclique number, which can be computed by the previous algorithm.

We first recall some notions of hypergraph \cite{Cooper12, Qi14}. A hypergraph means an undirected simple $m$-uniform hypergraph $G=(V,E)$ with
vertex set $V=\{1,2,\cdots,n\}$, and edge set $E = \{e_1,e_2,\cdots, e_k\}$ with $e_p \subseteq V$
for $p\in [k]$. By $m$-uniformity, we mean that for every edge $e\in E$, the cardinality $|e|$ of
$e$ is equal to $m$. A 2-uniform hypergraph is typically called graph. Throughout this paper,
we focus on $m\geq 3$ and $n \geq m$. Moreover, since the trivial hypergraph (i.e., $E=\emptyset$) is of
less interest, we consider only hypergraphs having at least one edge (i.e., nontrivial) in this
section.
\bd\label{def51}{\bf (Coclique number of a hypergraph)}
The coclique of an $m$-uniform hypergraph $G$ is a set of
vertices such that any of its $m$ vertex subset is not an edge of $G$, and the largest cardinality
of a coclique of $G$ is called the coclique number of $G$, denoted by $\omega(G)$.
\ed

By Definition \ref{def51}, we can easily get the following results.
\bp\label{prop51} Suppose $G=(V,E)$ is a nontrivial $m$-uniform hypergraph. Let $|V|=n$. Then,
the coclique number $\omega(G)$ of $G$ satisfies that $m-1\leq \omega(G) \leq n-1$.
\ep

The following definition for the adjacency tensor was proposed by Cooper and Dutle \cite{Cooper12}, which is important in the following
analysis.
\bd\label{def52}{\bf (Adjacency tensor of a hypergraph)} Let $G=(V,E)$ be an $m$-uniform hypergraph where $V=\{1,2,\cdots,n\}$. The
adjacency tensor of $G$ is defined as the $m$-th order $n$ dimensional tensor $\mathcal{A}$ with
$$
a_{i_1i_2\cdots i_m}=\left\{
\begin{array}{cll}
\frac{1}{(m-1)!} & \{i_1,i_2,\cdots,i_m\}\in E, \\
0 & {\bf otherwise}.
\end{array}
\right.
$$
\ed

\bt\label{them51} Let $G=(V,E)$ be an $m$-uniform hypergraph. Suppose $|V|=n$ and $G$ is nontrivial.
Let $\omega(G)$ denote the coclique number of $G$. Then $\omega(G)^{m-1}$ is equal to the optimal
value of the following problem:
$$
\begin{array}{clll}
{\rm (P)}~~~& \max & \langle \mathcal{X}, \mathcal{E}\rangle \\
& {\rm s.t.} & \mathcal{X}_{i_1i_2\cdots i_m}=0,~~\{i_1,i_2,\cdots,i_m\}\in E, \\
&    & \langle \mathcal{X}, \mathcal{I}\rangle=1, \\
&   & \mathcal{X}\in \mathbb{CP}_{m,n},
\end{array}
$$
where $\mathcal{E}$ is a all one tensor with order $m$ and dimension $n$.
\et
\proof By the fact that $\mathbb{CP}_{m,n}$ is a convex cone \cite{QXX14}, it is apparent that the feasible set of problem
(P) is also convex. So its optimal value will be attained at an extreme point, i.e., there is ${\bf x}^*\in \mathbb{R}^n_+$
such that $f^*=\langle ({\bf x}^*)^m, \mathcal{E} \rangle$, where $f^*$ is the optimal value of (P). Constraint conditions of
(P) implies that $\|{\bf x}\|_m=1$ and the support set $S^*$ of ${\bf x}^*$ is a coclique of $G$. By the optimal conditions of
(P), we can easily get that all nonzero entries of ${\bf x}^*$ are equal, which means that for any $i\in [n]$
$$
{\bf x}^*_i=
\left\{
\begin{array}{cll}
\frac{1}{\sqrt[m]{|S^*|}} & i\in S^*, \\
0 & {\rm otherwise}.
\end{array}
\right.
$$
Thus, the optimal value $f^*=\langle ({\bf x}^*)^m, \mathcal{E} \rangle=({\bf e}^T{\bf x}^*)^m=|S^*|^{m-1}$, which
implies that $S^*$ must be the maximum coclique and the desired result holds.  \qed

\bt\label{them52}
Assume $G$ and $\omega(G)$ are defined as in Theorem \ref{them51}. Let $\mathcal{A}$ be the adjacency tensor of hypergraph $G$.
Then, it holds that
$$
\omega(G)^{m-1}\leq \min_{\lambda \in \mathbb{N}}\{\lambda~|~\lambda(\mathcal{A}+\mathcal{I})-\mathcal{E}\in \mathbb{COP}_{m,n}\}.
$$
\et
\proof Since $\mathcal{X}\in \mathbb{CP}_{m,n}$ is nonnegative, by the proof of Theorem \ref{them51}, we know that problem (P) is equivalent to
$$
\begin{array}{clll}
 \max & \langle \mathcal{X}, \mathcal{E}\rangle \\
 {\rm s.t.} & \langle\mathcal{X},\mathcal{A}\rangle=0, \\
  & \langle \mathcal{X}, \mathcal{I}\rangle=1, \\
  & \mathcal{X}\in \mathbb{CP}_{m,n},
\end{array}
$$
which can be relaxed to the problem such that
$$
\begin{array}{clll}
{\rm (P')} & \max & \langle \mathcal{X}, \mathcal{E}\rangle \\
& {\rm s.t.} & \langle\mathcal{X},\mathcal{A}+\mathcal{I}\rangle=1, \\
&  & \mathcal{X}\in \mathbb{CP}_{m,n}.
\end{array}
$$
Then, the dual problem of ${\rm (P')}$ is that
$$
\min_{\lambda \in \mathbb{N}}\{\lambda~|~\lambda(\mathcal{A}+\mathcal{I})-\mathcal{E}\in \mathbb{COP}_{m,n}\}.
$$
From the well known weak duality theorem, we have
$$
\begin{aligned}
\omega(G)^{m-1}\leq & \max \{\langle \mathcal{X}, \mathcal{E}\rangle ~|~\langle\mathcal{X},\mathcal{A}+\mathcal{I}\rangle=1, \mathcal{X}\in \mathbb{CP}_{m,n} \} \\
\leq & \min_{\lambda \in \mathbb{N}}\{\lambda~|~\lambda(\mathcal{A}+\mathcal{I})-\mathcal{E}\in \mathbb{COP}_{m,n}\},
\end{aligned}
$$
which implies the desired results hold. \qed
By Proposition \ref{prop51} and Theorem \ref{them52}, we can try finitely many iterations to get an upper bound for the coclique number of a given uniform hypergraph by Algorithms 2. For example, for an $m$-uniform hypergraph $G=(V,E)$ with $V=[n]$, if there is $k\in [n]$ such that
$$
k^{m-1}(\mathcal{A}+\mathcal{I})-\mathcal{E}\in \mathbb{COP}_{m,n},~~(k-1)^{m-1}(\mathcal{A}+\mathcal{I})-\mathcal{E}\notin \mathbb{COP}_{m,n},
$$
then we know that the coclique number of $G$ satisfies $\omega(G)\leq k$.

\setcounter{equation}{0}
\section{Numerical results}

In this section, we report some preliminary numerical results of Algorithm 2, where the subcone $\mathbb{M}$ is chosen according to Section 4.2; and we use YALMIP \cite{Yalmip1,Yalmip2} and Sedumi \cite{sturm99} to solve the resulted semidefinite programs. All experiments are finished in Matlab2014a on a HP Z800 Workstation with Intel(R) Xeno(R) CPU X5680 @ 3.33GHz 3.33 GHz and 48 GB of RAM. All experiments are divided into the following two parts.

{\bf Part 1 (Copositivity detection)}. In this part, we implement Algorithm 5.2 to detect whether a tensor is copositive or not by using several examples tested in \cite{chen16}.
\begin{example}\label{example1}
We test the tensor $\mathcal{A}$ in the following form:
\begin{eqnarray}\label{E-exam1-1}
\mathcal{A}=\eta \mathcal{I}-\mathcal{B},
\end{eqnarray}
(i) Suppose that $\mathcal{A}\in S_{3,3}$ (or $\mathcal{A}\in S_{4,4}$) is given by (\ref{E-exam1-1}), where $\mathcal{B}\in S_{3,3}$ (or $\mathcal{B}\in S_{4,4}$) is a tensor of ones and $\eta$ is specified in the table of our numerical results.\\
(ii) Suppose that $\mathcal{A}\in S_{m,n}$ is given by (\ref{E-exam1-1}), where $\mathcal{B}\in S_{m,n}$ is randomly generated with all its elements are in the interval $(0,1)$.
\end{example}

The numerical results of testing the tensor $\mathcal{A}$ defined by Example \ref{example1}(i) are given in Table 1, where ``$\rho$" denotes the spectral radius of the tensor $\mathcal{B}$, ``IT" denotes the number of iterations, ``CPU(s)" denotes the CPU time in seconds, and ``Result" denotes the output result in which ``No" denotes the output result that the tested tensor is not copositive and ``Yes" denotes the output result that the tested tensor is copositive.

\begin{table}[ht] 
  \caption{The numerical results of the problem in Example \ref{example1}(i)}
  \begin{center}
    \begin{tabular}[c]
      {| c | c | c | c | c | c | c|}
      \hline
      $m$  & $n$   &$\rho$ & $\eta$ & IT      & CPU(s) & Result\\
      \hline
           &       &       & 1      & 2       &3.37     & No   \\  \cline{4-7}
           &       &       & 8.99   & 20      &81.2     & No   \\ \cline{4-7}
        3  & 3     & 9     & 9      & $>100$  &          &      \\  \cline{4-7}
           &       &       & 9.01   & 1       &0.92     & Yes   \\  \cline{4-7}
           &       &       & 19     & 1       &0.967    & Yes   \\ \hline
           &       &       & 10     & 8       &32.6     & No   \\ \cline{4-7}
        4  & 4     & 64    & 64     & $>100$  &         &    \\ \cline{4-7}
           &       &       & 74     & 1       &1.33     & Yes   \\ \hline
    \end{tabular}
  \end{center}
\end{table}

The numerical results of testing the tensor $\mathcal{A}$ defined by Example \ref{example1}(ii) are given in Table 2, where the spectral radius $\rho$ of every tensor $\mathcal{B}$ is computed by the higher order power method. In our experiments, for the same $m$ and $n$, we generate randomly every tested problem 10 times. In Table 2, ``MinIT" and ``MaxIT" denote the minimal number and the maximal number of iterations among ten times experiments for every tested problem, respectively, ``MinCPU(s)" and ``MaxCPU(s)" denote the smallest and the largest CPU times in second among ten times experiments for every tested problem, respectively, ``Nyes" denotes the number of the output results that the tested tensors are copositive, and ``Nno" denotes the number of the output results that the tested tensors are not copositive.

\begin{table}[ht] 
  \caption{The numerical results of the problem in Example \ref{example1}(ii)}
  \begin{center}
    \begin{tabular}[c]
      {| c | c | c | c | c | c | c| c | c |}
      \hline
      $m$  & $n$   & $\eta$       & MinIT   & MaxIT       &MinCPU(s) &MaxCPU(s) & Nyes  & Nno\\
      \hline
        3  & 3     & $\rho-1$  & 8       & 8          &18.7825    &30.9506     &       & 10    \\  \cline{3-9}
           &       & $\rho+1$  & 1       & 1          &0.98       &2.96        & 10    &       \\  \cline{3-9}
           &       & $\rho+10$ & 1       & 1          &1.02       &3.27        & 10    &      \\  \hline
        3  & 4     & $\rho-1$  & 16      & 16         &70.9337    &82.2125     &       & 10    \\  \cline{3-9}
           &       & $\rho+1$  & 1       & 1          &1.95       &3.8844      & 10    &       \\  \cline{3-9}
           &       & $\rho+10$ & 1       & 1          &2.184      &4.1808      & 10    &      \\  \hline
        4  & 3     & $\rho-1$  & 16      & 22         &72.5249    &106.3147    &       & 10    \\  \cline{3-9}
           &       & $\rho+1$  & 1       & 1          &1.5912     &3.2916      & 10    &       \\  \cline{3-9}
           &       & $\rho+10$ & 1       & 1          &1.7472     &3.354       & 10    &      \\  \hline
        4  & 4     & $\rho-1$  & 11      & 11         &47.2683    &52.3695     &       & 10    \\  \cline{3-9}
           &       & $\rho+1$  & 1       & 1          &1.7628     &3.8688      & 10    &       \\  \cline{3-9}
           &       & $\rho+10$ & 1       & 1          &1.8252     &3.978       & 10    &      \\  \hline
        6  & 3     & $\rho-1$  & 17      & 27         &97.6254    &157.3114    &       & 10    \\  \cline{3-9}
           &       & $\rho+1$  & 1       & 1          &2.574      &2.9016      & 10    &       \\  \cline{3-9}
           &       & $\rho+10$ & 1       & 1          &2.5116     &2.6676      & 10    &      \\  \hline
    \end{tabular}
  \end{center}
\end{table}

From Tables 1 and 2, it is easy to see that the concerned tensors can be correctly tested with few number of iterations. In particular, for every strictly copositive tensor we tested, only one step iteration reaches the right conclusion.

\begin{example}\label{example2}
We test the following three tensors:

(i) Suppose that $\mathcal{A}\in S_{6,3}$ is given by
\begin{eqnarray*}
\left\{\begin{array}{ll}
\sum_{i_1i_2i_3i_4i_5i_6\in S_{\pi(111122)}}a_{i_1i_2i_3i_4i_5i_6}=1,\\
\sum_{i_1i_2i_3i_4i_5i_6\in S_{\pi(112222)}}a_{i_1i_2i_3i_4i_5i_6}=1,\\
a_{333333}=1,\\
\sum_{i_1i_2i_3i_4i_5i_6\in S_{\pi(112233)}}a_{i_1i_2i_3i_4i_5i_6}=-3;
\end{array}\right.
\end{eqnarray*}
(ii) Suppose that $\mathcal{A}\in S_{6,3}$ is given by
\begin{eqnarray*}
\left\{\begin{array}{ll}
a_{111111}=1,\;\; a_{222222}=1,\;\; a_{333333}=1,\\
\sum_{i_1i_2i_3i_4i_5i_6\in S_{\pi(111122)}}a_{i_1i_2i_3i_4i_5i_6}=-1,\\
\sum_{i_1i_2i_3i_4i_5i_6\in S_{\pi(112222)}}a_{i_1i_2i_3i_4i_5i_6}=-1,\\
\sum_{i_1i_2i_3i_4i_5i_6\in S_{\pi(111133)}}a_{i_1i_2i_3i_4i_5i_6}=-1,\\
\sum_{i_1i_2i_3i_4i_5i_6\in S_{\pi(113333)}}a_{i_1i_2i_3i_4i_5i_6}=-1,\\
\sum_{i_1i_2i_3i_4i_5i_6\in S_{\pi(222233)}}a_{i_1i_2i_3i_4i_5i_6}=-1,\\
\sum_{i_1i_2i_3i_4i_5i_6\in S_{\pi(223333)}}a_{i_1i_2i_3i_4i_5i_6}=-1,\\
\sum_{i_1i_2i_3i_4i_5i_6\in S_{\pi(112233)}}a_{i_1i_2i_3i_4i_5i_6}=3;
\end{array}\right.
\end{eqnarray*}
(iii) Suppose that $\mathcal{A}\in S_{6,3}$ is given by
\begin{eqnarray*}
\left\{\begin{array}{ll}
\sum_{i_1i_2i_3i_4i_5i_6\in S_{\pi(111122)}}a_{i_1i_2i_3i_4i_5i_6}=1,\\
\sum_{i_1i_2i_3i_4i_5i_6\in S_{\pi(222233)}}a_{i_1i_2i_3i_4i_5i_6}=1,\\
\sum_{i_1i_2i_3i_4i_5i_6\in S_{\pi(333311)}}a_{i_1i_2i_3i_4i_5i_6}=1,\\
\sum_{i_1i_2i_3i_4i_5i_6\in S_{\pi(112233)}}a_{i_1i_2i_3i_4i_5i_6}=-3.
\end{array}\right.
\end{eqnarray*}
\end{example}

The corresponding polynomials of the above tensors are famous Motzkin polynomial, Robinson polynomial and Choi-Lam polynomial, respectively. It is easy to see that the above three tensors are copositive, but not strictly copositive. We use Algorithm 5.2 to test the tensor $\mathcal{A}+\sigma \mathcal{E}$ with $\sigma>0$, and the numerical results are listed in Table 3.

\begin{table}[ht] 
  \caption{The numerical results of the problem in Example \ref{example2}}
  \begin{center}
    \begin{tabular}[c]
      {| c | c | c | c |}
      \hline
                & Example 6.3(i)   & Example 6.3(ii)     & Example 6.3(iii)   \\  \cline{2-4}
      $\sigma$  & IT/CPU(s)        & IT/CPU(s)           & IT/CPU(s)           \\
      \hline
        0.01    & 3/14.3           & 11/60               &   5/27.4     \\  \hline
        0.001   & 19/92.8          & 27/152              &   17/99.5     \\  \hline
        0.0001  & 55/291           & 67/379              &   35/209     \\  \hline
    \end{tabular}
  \end{center}
\end{table}

All tensors tested in Examples \ref{example1} and \ref{example2} were tested in \cite{chen16}. Compared the numerical results shown in Tables 1-3 with those given in \cite{chen16}, choosing
$\mathbb{M}= \mathbb{H}$ requires the least number of iterations but each iteration is so costly that the
overall runtime is in most cases still higher than $\mathbb{M}=\mathbb{N}^+_{m,n}$.

{\bf Part 2 (Illustration of Theorem 5.2)}. As said in the last section, for an $m$-uniform hypergraph $G=(V,E)$ with $V=[n]$, if there is $k\in [n]$ such that
$$
k^{m-1}(\mathcal{A}+\mathcal{I})-\mathcal{E}\in \mathbb{COP}_{m,n}~~{\bf or}~~(k-1)^{m-1}(\mathcal{A}+\mathcal{I})-\mathcal{E}\notin \mathbb{COP}_{m,n},
$$
then we know that the coclique number of $G$ satisfies $\omega(G)\leq k$. By this way, we can compute the coclique number of an $m$-uniform hypergraph. Conversely, if the coclique number of an $m$-uniform hypergraph is known, we can also check the main result obtained in Section 5. In this part, we illustrate Theorem 5.2 by constructing two examples.

\begin{example}\label{example-1}
Let $V=\{1,2,3\}$ and $E$ be a set of subsets of $V$. Let $G=(V,E)$ be a $3$-uniform
hypergraph. If $V_0=\{1\}$, $V_1=\{2,3\}$ and $E=\left\{1,2,3\right\}$, then $G$ is a hyper-star.
\end{example}

The adjacency tensor of $G$ is as follows:
\begin{eqnarray*}
A(:,:,1) = \left(\begin{array}{ccc}
0 & 0 & 0\\
0 & 0 & \frac{1}{2}\\
0 & \frac{1}{2} & 0
\end{array}\right),
A(:,:,2) = \left(\begin{array}{ccc}
0 & 0 & \frac{1}{2}\\
0 & 0 & 0\\
\frac{1}{2} & 0 & 0
\end{array}\right),
A(:,:,3) = \left(\begin{array}{ccc}
0 & \frac{1}{2} & 0\\
\frac{1}{2} & 0 & 0\\
0 & 0 & 0
\end{array}\right).
\end{eqnarray*}
By Algorithm 2, we can obtain that $4(\mathcal{A}+\mathcal{I})-\mathcal{E}\notin \mathbb{COP}_{m,n}$ (This can be also seen by $f({\bf x})=4(\mathcal{A}+\mathcal{I}){\bf x}^3-\mathcal{E}{\bf x}^3=-3<0$ when ${\bf x}=(1,1,1)\in \mathbb{R}^3$). So, from the monotonicity of $\lambda(\mathcal{A}+\mathcal{I})-\mathcal{E}$, we know that
$\min_{\lambda \in \mathbb{N}}\{\lambda~|~\lambda(\mathcal{A}+\mathcal{I})-\mathcal{E}\in \mathbb{COP}_{m,n}\}>4$. By the fact that $\omega(G)=n-1=2$, we have
$$
\omega(G)^{m-1}\leq \min_{\lambda \in \mathbb{N}}\{\lambda~|~\lambda(\mathcal{A}+\mathcal{I})-\mathcal{E}\in \mathbb{COP}_{m,n}\},
$$
which verifies Theorem 5.2.

\begin{example}\label{example-2}
Let $V=\{1,2,3,4\}$ and $E$ be a set of subsets of $V$. Let $G=(V,E)$ be a $4$-uniform
hypergraph. If $V_0=\{1\}$, $V_1=\{2,3,4\}$ and $E=\left\{\{1,2,3,4\}\right\}$, then $G$ is a hyper-star.
\end{example}

The coefficients of the adjacency tensor $\mathcal A$ of $G$ are as follows:
\begin{eqnarray*}
&&A(:,:,1,1) = \left(\begin{array}{cccc}
0 & 0 & 0 & 0 \\
0 & 0 & 0 & 0\\
0 & 0 & 0 & 0 \\
0 & 0 & 0 & 0
\end{array}\right),
A(:,:,2,1) = \left(\begin{array}{cccc}
0 & 0 & 0 & 0 \\
0 & 0 & 0 & 0 \\
0 & 0 & 0 & \frac{1}{6} \\
0 & 0 & \frac{1}{6} & 0
\end{array}\right),
A(:,:,3,1) = \left(\begin{array}{cccc}
0 & 0 & 0 & 0 \\
0 & 0 & 0 & \frac{1}{6} \\
0 & 0 & 0 & 0 \\
0 & \frac{1}{6} & 0 & 0
\end{array}\right),\\
&&A(:,:,4,1) = \left(\begin{array}{cccc}
0 & 0 & 0 & 0 \\
0 & 0 & \frac{1}{6} & 0\\
0 & \frac{1}{6} & 0 & 0\\
0 & 0 & 0 & 0
\end{array}\right),
A(:,:,1,2)  = \left(\begin{array}{cccc}
0 & 0 & 0 & 0 \\
0 & 0 & 0 & 0 \\
0 & 0 & 0 & \frac{1}{6} \\
0 & 0 & \frac{1}{6} & 0
\end{array}\right),
A(:,:,2,2) = \left(\begin{array}{cccc}
0 & 0 & 0 & 0 \\
0 & 0 & 0 & 0\\
0 & 0 & 0 & 0 \\
0 & 0 & 0 & 0
\end{array}\right),\\
&&A(:,:,3,2)  = \left(\begin{array}{cccc}
0 & 0 & 0 & \frac{1}{6} \\
0 & 0 & 0 & 0 \\
0 & 0 & 0 & 0 \\
\frac{1}{6} & 0 & 0 & 0
\end{array}\right),
A(:,:,4,2)  = \left(\begin{array}{cccc}
0 & 0 & \frac{1}{6} & 0 \\
0 & 0 & 0 & 0 \\
\frac{1}{6} & 0 & 0 & 0\\
0 & 0 & 0 & 0
\end{array}\right),
A(:,:,1,3)  = \left(\begin{array}{cccc}
0 & 0 & 0 & 0 \\
0 & 0 & 0 & \frac{1}{6} \\
0 & 0 & 0 & 0 \\
0 &\frac{1}{6} & 0 & 0
\end{array}\right),\\
&&A(:,:,2,3)  = \left(\begin{array}{cccc}
0 & 0 & 0 & \frac{1}{6} \\
0 & 0 & 0 & 0 \\
0 & 0 & 0 & 0 \\
\frac{1}{6} & 0 & 0 & 0
\end{array}\right),
A(:,:,3,3) = \left(\begin{array}{cccc}
0 & 0 & 0 & 0 \\
0 & 0 & 0 & 0\\
0 & 0 & 0 & 0 \\
0 & 0 & 0 & 0
\end{array}\right),
A(:,:,4,3)  = \left(\begin{array}{cccc}
0 & \frac{1}{6} & 0 & 0\\
\frac{1}{6} & 0 & 0 & 0\\
0 & 0 & 0 & 0 \\
0 & 0 & 0 & 0
\end{array}\right),\\
&&A(:,:,1,4) = \left(\begin{array}{cccc}
0 & 0 & 0 & 0 \\
0 & 0 & \frac{1}{6} & 0\\
0 & \frac{1}{6} & 0 & 0\\
0 & 0 & 0 & 0
\end{array}\right),
A(:,:,2,4)  = \left(\begin{array}{cccc}
0 & 0 & \frac{1}{6} & 0 \\
0 & 0 & 0 & 0 \\
\frac{1}{6} & 0 & 0 & 0\\
0 & 0 & 0 & 0
\end{array}\right),
A(:,:,3,4)  = \left(\begin{array}{cccc}
0 & \frac{1}{6} & 0 & 0\\
\frac{1}{6} & 0 & 0 & 0\\
0 & 0 & 0 & 0 \\
0 & 0 & 0 & 0
\end{array}\right),\\
&&A(:,:,4,4) = \left(\begin{array}{cccc}
0 & 0 & 0 & 0 \\
0 & 0 & 0 & 0\\
0 & 0 & 0 & 0 \\
0 & 0 & 0 & 0
\end{array}\right),
\end{eqnarray*}
By running program of Algorithm 2, we obtain $8(\mathcal{A}+\mathcal{I})-\mathcal{E} \notin \mathbb{COP}_{m,n}$, which implies that
$$\min_{\lambda \in \mathbb{N}}\{\lambda~|~\lambda(\mathcal{A}+\mathcal{I})-\mathcal{E}\in \mathbb{COP}_{m,n}\}\geq27.$$
Since $\omega(G)=3$, it holds that
$$
\omega(G)^{m-1}\leq \min_{\lambda \in \mathbb{N}}\{\lambda~|~\lambda(\mathcal{A}+\mathcal{I})-\mathcal{E}\in \mathbb{COP}_{m,n}\},
$$
and hence, Theorem 5.2 is true.

\setcounter{equation}{0}
\section{Checking vacuum stability for $\mathbb{Z}_3$ scalar dark matter}
Kannike \cite{KK16} studied the vacuum stability of a general scalar potential of a few fields,
and explicit vacuum stability conditions for more complicated potentials are given. In \cite{KK16}, one important
physical example is given by scalar dark matter stable under $\mathbb{Z}_3$ discrete group. The most general
scalar quartic potential of the standard model(SM) Higgs $H_1$, an inert doublet $H_2$ and a complex singlet
$S$ which is symmetric under a $\mathbb{Z}_3$ group is
\begin{equation}\label{e71}
\begin{aligned}
V(h_1,h_2,S)=& \lambda_1|H_1|^4+\lambda_2|H_2|^4+\lambda_3|H_1|^2|H_2|^2+\lambda_4(H_1^{\dag}H_2)(H_2^{\dag}H_1)+\lambda_S|S|^4
+\lambda_{S1}|S|^2|H_1|^2\\
&+\lambda_{S2}|S|^2|H_2|^2+\frac{1}{2}(\lambda_{S12}S^2H_1^{\dag}H_2+\lambda_{S12}^*{S^{\dag}}^2H_2^{\dag}H_1) \\
=&\lambda_1h_1^4+\lambda_2h_2^4+\lambda_3h_1^2h_2^2+\lambda_4\rho^2h_1^2h_2^2+\lambda_Ss^4+\lambda_{S1}s^2h_1^2
+\lambda_{S2}s^2h_2^2-|\lambda_{S12}|\rho s^2h_1h_2 \\
\equiv& \lambda_Ss^4+M^2(h_1,h_2)s^2+V(h_1,h_2),
\end{aligned}
\end{equation}
where $M^2(h_1,h_2):=\lambda_{S1}s^2h_1^2 +\lambda_{S2}s^2h_2^2-|\lambda_{S12}|\rho s^2h_1h_2$ and $V(h_1,h_2):=V(h_1,h_2,0)$.
Here, in physical sense, the variables $h_1, h_2$ and $s$ are nonnegative since they are magnitudes of scalar fields, so the coupling tensor $\mathcal V$ of coefficients of
(\ref{e71}) has to be copositive. This has to hold for all values of the extra parameter $\rho$ ranges from 0 to 1, so the potential has to be minimized or scanned over it. Now, we give the explicit form for the coupling tensor of (\ref{e71}) as ${\mathcal V}=(V_{i_1i_2i_3i_4})$, which is an order 4 dimension 3 real symmetric tensor:
$$
V_{1111}=\lambda_1,~~V_{2222}=\lambda_2,~~V_{3333}=\lambda_S
$$
$$
V_{1122}=\frac{1}{6}(\lambda_3+\lambda_4\rho^2), V_{1133}=\frac{1}{6}\lambda_{S1}, V_{2233}=\frac{1}{6}\lambda_{S2}, V_{1233}=-\frac{1}{12}|\lambda_{S12}|
$$
and $V_{i_1i_2i_3i_4}=0$ for the others. Then, by Algorithm 2, we give a series of explicit coefficients and check the vacuum stability of the potential (\ref{e71}).

As to $\lambda$'s in the entries of $\mathcal V$, in particle physics all calculated quantities are expanded in series of $\lambda_i/(4 \pi)$. Due to the perturbativity requirement of these series, the absolute values of the $\lambda$ coefficients must be no larger than $4 \pi$. On the other hand, for the coupling tensor to be copositive, the diagonal entries have to be nonnegative.   Hence, we can take from the beginning that $0 \leq V_{1111}, V_{2222}, V_{3333} \leq 4 \pi$.

Then, because the rest of the entries of $\mathcal V$ are a $\lambda$ paremeter times some coefficients, their lower and upper bounds should be accordingly changed. So $-2 \times 4 \pi/6 \leq V_{1122} \leq 2 \times 4 \pi/6$ (with an extra factor $2$ because it is the sum of two $\lambda$'s), $-4 \pi/6 \leq V_{1133} \leq 4 \pi/6$, $-4 \pi/6 \leq V_{1133} \leq 4 \pi/6$, $-4 \pi/6 \leq V_{2233} \leq 4 \pi/6$, and $-4 \pi/12 \leq V_{1233} \leq 0$.

When $\rho\neq 0$, Kannike \cite{KK16} obtained that the conditions for the potential (\ref{e71}) symmetric under a $\mathbb{Z}_3$ to be bounded from below are
\begin{eqnarray}\label{e72}
\left\{\begin{array}{l}
\lambda_S>0,\\
V(h_1,h_2)>0,\\
0<h_1^2<1, 0<h_2^2<1, 0<s^2<1,\;\mbox{\rm and}\; 0<\rho^2<1\quad \Longrightarrow \quad V_{\min}>0,
\end{array}\right.
\end{eqnarray}
where
\begin{eqnarray}\label{e73}
\begin{array}{rcl}
\rho&=&\left(|\lambda_{S12}|s^2\right)/\left(2\lambda_4h_1h_2\right),\\
h_1^2&=&\frac 12\left\{(2\lambda_2-\lambda_3)(4\lambda_S\lambda_4-|\lambda_{S12}|^2)+2\lambda_4\left[
(\lambda_3+\lambda_{S1})\lambda_{S2}-2\lambda_2\lambda_{S1}-\lambda_{S2}^2\right]\right\}/t,\\
h_2^2&=&\frac 12\left\{(2\lambda_1-\lambda_3)(4\lambda_S\lambda_4-|\lambda_{S12}|^2)+2\lambda_4\left[
(\lambda_3+\lambda_{S2})\lambda_{S1}-2\lambda_1\lambda_{S2}-\lambda_{S1}^2\right]\right\}/t,\\
s^2&=& \lambda_4\left(4\lambda_1\lambda_2-\lambda_3^2-2\lambda_1\lambda_{S2}-2\lambda_2\lambda_{S1}+\lambda_3(\lambda_{S1}+\lambda_{S2})\right)/t,\\
V_{\min}&=& \frac{1}{4}\left[(4\lambda_1\lambda_2-\lambda_3^2)(4\lambda_S\lambda_4-|\lambda_{S12}|^2)-4\lambda_4(\lambda_1\lambda_{S2}^2 +\lambda_2\lambda_{S1}^2-\lambda_3\lambda_{S1}\lambda_{S2})\right]/t
\end{array}
\end{eqnarray}
with
\begin{eqnarray*}
\begin{array}{rcl}
t&:=&(\lambda_1+\lambda_2-\lambda_3)\times (4\lambda_S\lambda_4-|\lambda_{S12}|^2)\\
&&+\lambda_4\left[4\lambda_1\lambda_2-\lambda_3^2-4\lambda_1\lambda_{S2}-4\lambda_2\lambda_{S1}
+2\lambda_3(\lambda_{S1}+\lambda_{S2})-(\lambda_{S1}-\lambda_{S2})^2\right].
\end{array}
\end{eqnarray*}
And the third formula in (\ref{e72}) is replaced by $V_{\rho=0}>0$ when $\rho=0$; and by $V_{\rho=1}>0$ when $\rho=1$.

Now, we implement Algorithm 2 to test the copositivity of the tensor defined by the potential (\ref{e71}), and the numerical results are listed in Table 4, where the values of $h_1^2, h_2^2, s^2, \rho$ and $V_{\min}$ are computed by (\ref{e73}), ``IT" denotes the number of iteration, ``CPU(s)" denotes the CPU time in seconds, ``Yes" denotes the output result that the tested tensor is copositive and ``No" denotes the output result that the tested tensor is not copositive.

\begin{table}
  \caption{The numerical results for the stability of the potential (\ref{e71})}
  \begin{center}
    \begin{tabular}[c]
      {| c | c | c | c | c | c | c| c | c | c | c | c| c | c | c | c |}
      \hline
    $\lambda_1$  & $\lambda_2$   & $\lambda_S$    & $\lambda_3$   & $\lambda_4$       & $\lambda_{S1}$ & $\lambda_{S2}$ & $\lambda_{S12}$ &
    $h_1^2$ & $h_2^2$ & $s^2$ & $\rho$ & $V_{\min}$ & IT  & CPU(s) & Result  \\
      \hline
 $\pi$  &  $\pi$   & $\pi$   & $\pi$   & $\pi$   &$\pi$   &$\pi$   & $\pi$  & 0.25 & 0.25 & 0.5  & 1    & 1.96 & 3  & 0.09 & Yes  \\  \hline
 $\pi$  &  $\pi$   & $\pi$   & $\pi$   & $\pi$   &-$\pi$  &-$\pi$  & 0      & 0.27 & 0.27 & 0.45 & 0    & 0.57 & 19 & 0.37 & Yes   \\  \hline
 $\pi$/4&  $\pi$/4 & $\pi$/4 & $\pi$   & $\pi$   &0       &0       & 4$\pi$ & 0.56 & 0.56 &-0.11 & 0.4  & 1.31 & 9  & 0.17 & No  \\  \hline
 $\pi$  &  $\pi$   & $\pi$/2 & $\pi$   & $\pi$   &0       &0       & 4$\pi$ & 0.64 & 0.64 &-0.27 & 0.86 & 3.00 & 5  & 0.08 & No \\ \hline
 $\pi$  &  $\pi$   & -$\pi$  & $\pi$   & $\pi$   &$\pi$   &$\pi$   & 4$\pi$ & 0.52 & 0.52 &-0.05 & 0.18 & 2.39 & 1  & 0.02 & No  \\  \hline
 $-\pi$ &  $\pi$   & $\pi$   & $\pi$   & $\pi$   &$\pi$   &$\pi$   & 4$\pi$ &-0.64 & 1.91 &-0.27 &0.49i & 4.57 & 1  & 0.02 & No \\  \hline
 $\pi$  &  $-\pi$  & $\pi$   & $\pi$   & $\pi$   &$\pi$   &$\pi$   & 4$\pi$ & 1.91 &-0.64 &-0.27 &0.49i & 4.57 & 1  & 0.02 & No \\ \hline
 $\pi$  &  $\pi$   & -$\pi$  & $\pi$   & $\pi$   &$\pi$   &$\pi$   & 3$\pi$ & 0.54 & 0.54 &-0.07 & 0.2  & 2.41 & 1  & 0.02 & No  \\ \hline
 $-\pi$ &  $\pi$   & $\pi$   & $\pi$   & $\pi$   &$\pi$   &$\pi$   & 3$\pi$ &-0.88 & 2.63 &-0.75 &0.74i & 5.69 & 1  & 0.02 & No \\  \hline
 $\pi$  &  $-\pi$  & $\pi$   & $\pi$   & $\pi$   &$\pi$   &$\pi$   & 3$\pi$ & 2.63 &-0.88 &-0.75 &0.74i & 5.69 & 1  & 0.02 & No \\ \hline
 $\pi$/4&  $\pi$/4 & $\pi$/4 & $\pi$/4 & $\pi$/4 &$\pi$/2 &$\pi$/2 & 4$\pi$ & 0.50 &0.50  &0.004 & 0.06 & 0.59 & 3  & 0.05 & No  \\  \hline
 $\pi$  &  $\pi$   & $\pi$   & $\pi$   & $\pi$   &$\pi$   &$\pi$   & $\pi$/2& 0.32 & 0.32 & 0.36 & 0.29 & 2.07 & 3  & 0.06 & Yes  \\  \hline
 2$\pi$ &  $\pi$   & 2$\pi$  & $\pi$   & $\pi$   &$\pi$   &$\pi$   & $\pi$/2& 0.20 & 0.59 & 0.21 & 0.15 & 2.51 & 3  & 0.06 & Yes  \\  \hline
 $\pi$  &  $\pi$   & 2$\pi$  & $\pi$   & $\pi$   &$\pi$   & 0      & $\pi$/2& 0.24 & 0.50 & 0.26 & 0.19 & 1.95 & 7  & 0.13 & Yes  \\  \hline
 $\pi$  &  $\pi$   & 2$\pi$  & $\pi$   & $\pi$   & 0      &$\pi$   & $\pi$/2& 0.50 & 0.24 & 0.26 & 0.19 & 1.95 & 7  & 0.14 & Yes  \\  \hline
 $2\pi$ &  $\pi$   & 2$\pi$  & $\pi$   & $\pi$   & 0      &$\pi$   & $\pi$/2& 0.25 & 0.49 & 0.26 & 0.18 & 2.34 & 7  & 0.14 & Yes \\  \hline
 $\pi$  & $2\pi$   & 2$\pi$  & $\pi$   & $\pi$   & 0      & $\pi$  & $\pi$/2& 0.60 & 0.10 & 0.31 & 0.32 & 2.02 & 7  & 0.14 & Yes \\  \hline
 $2\pi$ & $2\pi$   & $\pi$   & $\pi$   & $\pi$   & 0      & $\pi$  & $\pi$/2& 0.29 & 0.08 & 0.62 & 0.99 & 1.97 & 7  & 0.14 & Yes \\  \hline
 $\pi$  & $\pi$    & $2\pi$  & 0       & $\pi$   & $\pi$  & 0      & $2\pi$ & 0.29 & 0.43 & 0.29 & 0.82 & 1.35 & 11 & 0.20 & Yes \\  \hline
 $\pi$/4& $\pi$/4  & $\pi$   & 0       & $\pi$   & $\pi$  & 0      & $2\pi$	& 0.29 & 0.57 & 0.14 & 0.35 & 0.45 & 7  & 0.14 & Yes \\  \hline
 $\pi$  & $\pi$    & $\pi$   & 0       & $\pi$   & -$\pi$ & 0      & $\pi$/2& 0.40 & 0.19 & 0.41 & 0.38 & 0.60 & 19 & 0.36 & Yes \\  \hline
 $\pi$  & $\pi$    & $\pi$   & 0       & $\pi$   & -$\pi$ & -$\pi$ & $2\pi$ & 0.17 & 0.17 &	0.67 & 4    &-0.52 & 17 & 0.31 & No \\  \hline
 $\pi$  & $\pi$    & $\pi$   & $\pi$   & $\pi$   & -$\pi$ & -$\pi$ & $2\pi$ & 0.14 & 0.14 & 0.71 & 5    &-0.45 & 17 & 0.31 & No \\  \hline
   \end{tabular}
  \end{center}
\end{table}

It should be noted that it is easy to check that the parameters are satisfied with the conditions (\ref{e72}) when our tested results are ``Yes", while the parameters are not satisfied with the conditions (\ref{e72}) when our tested results are ``No". We also tested other cases, the computation effect is similar. These imply that our algorithm is efficient and applicable to such physical problems.

\setcounter{equation}{0}
\section{Conclusions}

In this paper, an alternative form of a
previously given algorithm for copositivity of high order tensors is given, and  applications of the proposed algorithm to test copositivity of the coupling tensor in a vacuum stability model in particle physics, and to compute the coclique number of a uniform hypergraph are presented. Several new conditions for copositivity of tensors based on the representative matrix of a simplex are proved.
We see that for the performance of this algorithm the choice of the set $\mathbb{M}$ is crucial, and it is
observed that verifying copositivity of tensors is much harder than verifying non-copositivity.

However, some interesting questions still need to study in the future:

{\bf 1}. Are there any better choices for the set $\mathbb{M}$ in Algorithm 2?

{\bf 2}. How to update the proposed method to make it available for copositive tensors but not strictly copositive?

\bigskip

{\bf Acknowledgment}   We are thankful to Kristjan Kannike for the discussion on the vacuum stability model.

\end{document}